\input amstex

\def\b1{\text{\bf 1}}

\def\CA{{\Cal A}}

\def\CC{{\Cal C}}
\def\CD{{\Cal D}}

\def\CL{{\Cal L}}
\def\CM{{\Cal M}}

\def\CO{{\Cal O}}
\def\CP{{\Cal P}}

\def\CT{{\Cal T}}
\def\CV{{\Cal V}}

\def\gr{\text{gr}}

\def\End{\text{End}}
\def\Hom{\text{Hom}}
\def\Sym{\text{Sym}}

\def\#{\,\check{}}

\def\fg{{\frak g}}

\def\fl{{\frak{l}}}

\def\CK{\Cal K}

\def\id{\text{id}}

\def\Spec{\text{Spec}}

\def\limleft{\mathop{\vtop{\ialign{##\crcr
  \hfil\rm lim\hfil\crcr
  \noalign{\nointerlineskip}\leftarrowfill\crcr
  \noalign{\nointerlineskip}\crcr}}}}
\def\limright{\mathop{\vtop{\ialign{##\crcr
  \hfil\rm lim\hfil\crcr
  \noalign{\nointerlineskip}\rightarrowfill\crcr
  \noalign{\nointerlineskip}\crcr}}}}

\def\hra{\hookrightarrow}
\def\iso{\buildrel\sim\over\rightarrow}

\def\LC{\mathop{\vtop{\ialign{##\crcr
  \hfil\rm $\CC$\hfil\crcr
  \noalign{\nointerlineskip}\rightarrowfill\crcr
  \noalign{\nointerlineskip}\crcr}}}}

\def\LC^0{\mathop{\vtop{\ialign{##\crcr
  \hfil\rm $\CC^0$\hfil\crcr
  \noalign{\nointerlineskip}\rightarrowfill\crcr
  \noalign{\nointerlineskip}\crcr}}}}

\parskip=6pt

\documentstyle{amsppt}
\document
\magnification=1100
\NoBlackBoxes

\document

 \centerline{\bf  REMARKS ON TOPOLOGICAL  ALGEBRAS}

\bigskip

\centerline {A.~Beilinson}

\centerline {University of Chicago}
\medskip

\centerline{\it To Victor Ginzburg on his 50th birthday}

\bigskip

 The note   complements   ``topological" aspects of the chiral algebras story from \cite{BD}. In its first section (which has a whiff of \cite{G} in it) we show that the basic chiral algebra format (chiral operations, etc.) has a precise analog in the setting of topological linear algebra. This provides, in particular, a natural explanation of the passage from chiral to topological algebras from \cite{BD} 3.6.
The second section is a brief discussion, in the spirit of \cite{BD} 3.9, of topological algebras similar to rings of chiral differential operators.  We also  correct some  errors from \cite{BD}.
 
 The writing was prompted by a talk of D.~Gaitsgory at the end of 2004; I am grateful to him for helpful discussions.

\head
1.  Topological tensor products and topological algebras.
\endhead

{\bf 1.1. The tensor products.}
For us, as in \cite{BD} 3.6.1,   ``topological vector space" is a $k$-vector space equipped with a linear topology assumed  (unless stated explicitly otherwise) to be complete and separated. The category of topological vector spaces is denoted by $\CT\! op$. This is an additive Karoubian $k$-category. 

The category $\CT\! op$ is   quasi-abelian in the sense of   \cite{S}. In particular, it is
naturally an exact category: the admissible monomorphisms are closed embeddings, the admissible  epimorphisms are open surjections.

{\it Remark.} If the topology of $B\in\CT\! op$ admits a countable base, then every
 short exact sequence $0\to A \to B\to C \to 0$  is split exact, i.e., $B$ is isomorphic to $A\oplus C$. 
 \medskip

Let $\{ V_i \}_{i\in I}$ be a finite non-empty collection of topological vector spaces. Consider the tensor product $\otimes V_i$. This is an abstract vector space; it carries several natural linear topologies. Namely:

(a)\enspace The  {\it $*$ topology} is formed by all subspaces $Q$ of $ \otimes V_i$ 
which satisfy the following property:  for every subset $J\subset I$ and vector 
$v\in \mathop\otimes\limits_{i\in I\smallsetminus J}V_i $ 
there exist open subspaces $P_j \subset V_j$, $j\in J$, such that $Q\supset (\otimes P_j )\otimes v$. Let $\otimes^* V_i$ be the corresponding completion. Then for any topological space $F$ a continuous morphism $\otimes^* V_i \to F$ is the same as a continuous polylinear map $\times  V_i \to F$.

(b)\enspace 
The  {\it $!$ topology} has base formed by vector subspaces  $\mathop\sum\limits_{i\in I} P_i \otimes (\mathop\otimes\limits_{i' \in I\smallsetminus \{ i\} } V_{i'})$ where $P_i \subset V_i$ are any open subspaces. The corresponding completion is denoted by $\otimes^! V_i$.
Thus $\otimes^! V_i = \limleft \otimes (V_i /P_i )$.

 Notice that  \cite{BD} (in particular, in \cite{BD} 3.6.1) used the notation $\hat{\otimes}V_i$; the reason for  the  change of notation will become clear later.

(c)\enspace  Suppose we have a linear order $\tau$ on $I$, i.e., an identification $\{ 1,\ldots ,n\} @>{\sim}>> I$. It defines on $\otimes V_i$  the {\it $\tau$ topology} formed by all subspaces $Q$ which satisfy the following property: for every $a\in \{ 1,\ldots , n\}$ and vector $v\in V_{\tau (a+1)}\otimes \ldots \otimes V_{\tau (n)}$ there exists an open subspace $P_a \subset V_{\tau (a)}$ such that $Q\supset V_{\tau (1)} \otimes \ldots V_{\tau (a-1)}\otimes P_a \otimes v$. The corresponding completion is denoted by $\vec{\otimes}^\tau V_i = V_{\tau (1)} \vec{\otimes}\ldots \vec{\otimes} V_{\tau (n)}$.

\medskip

We refer to morphisms $\otimes V_i \to F$ which are continuous with respect to the ! topology, i.e., morphisms $\otimes^! V_i \to F$,  as $\otimes^!$-continuous polylinear maps. Same for the other tensor products.\footnote{So $\otimes^*$-continuous polylinear maps are the same as continuous polylinear maps.}

{\it Remarks.} (i) Let $U, V$ be topological vector spaces. Suppose first that  $V$ is discrete. Then $U\otimes^* V = U\vec{\otimes} V$ is equal to $U\otimes V$ equipped with the ind-topology. Precisely, write $V= \limright V_\alpha$ where $V_\alpha $ runs the directed set of finite-dimensional subspaces of $V$; then $U\otimes V = \limright U\otimes V_\alpha$, and each $U\otimes V_\alpha$ carries an evident topology (as the product of finitely many copies of $U$), which defines the said inductive limit topology on $U\otimes V$. 

If $V$ is arbitrary, then $U\vec{\otimes}V = \limleft U\vec{\otimes} (V/P)$, the projective limit is taken along the set of all open subspaces $P\subset V$.

(ii) Suppose we have $A\in \CT\!op$ and an associative bilinear product $\cdot : A\otimes A\to A$. Then $\cdot$ is $\otimes^*$-continuous  if and only if the product map $A\times A\to A$ is continuous; it is $\vec{\otimes}$-continuous  if and only if it is continuous and the open left ideals  form a base of the topology of $A$;\footnote{Or, equivalently,  the open left ideals  form a base of the topology  and for every $r\in A$ either of endomorphisms  $a \mapsto ar$ or $a\mapsto ar-ra$ of $A$ is continuous. } and it is $\otimes^!$-continuous  if and only if  the open two-sided ideals form a base of the topology.

\medskip

The tensor products $\otimes^*$ and $\otimes^!$ are  commutative and associative; they define tensor structures on $\CT\!op$  which we denote by $\CT\!op^!$ and $\CT\!op^*$. The tensor product $\vec{\otimes}$ is associative but {\it not} commutative; it defines a monoidal structure on $\CT\!op$.

{\it Exercise.} The tensor products are exact.

The $*$ topology on $\otimes V_i$ is stronger than each of the $\tau$ topologies, which in turn are stronger than the $!$ topology, so  the identity map for $\otimes V_i$ gives rise to the natural continuous morphisms  $$\otimes^* V_i \to \vec{\otimes}^\tau V_i \to \otimes^! V_i . \tag 1.1.1$$

\proclaim{\quad Lemma}  The $*$ topology is equal to the supremum of the $\tau$ topologies for all linear orders $\tau$ on $I$. The $!$ topology is the infimum of the $\tau$ topologies. I.e., the arrows in (1.1.1) are, respectively,  admissible mono- and  epimorphism.
  \hfill$\square$
\endproclaim

\proclaim{\quad Corollary} For any pair $U$, $V$ of topological vector spaces the short sequence $$ 0\to U \otimes^* V \to U \vec{\otimes} V \oplus V \vec{\otimes} U \to U\otimes^! V\to 0, \tag 1.1.2$$ where the left arrow comes from the diagonal map and the right one from the difference of the projections, is exact.  \hfill$\square$
\endproclaim

{\bf 1.2.} {\it Example.} Suppose  $F$, $G$ are Tate vector spaces (see e.g.~\cite{BD} 2.7.7, \cite{D} 3.1). Let $F^*$ be the dual Tate space to $F$, so $F^* \otimes G$ is the vector space of all continuous maps $F\to G$ with finite-dimensional image. Then $F^* \otimes^! G$ identifies naturally with the space $\Hom (F,G)$ of all continuous linear maps $F\to G$, $F^* \vec{\otimes} G$ with the  space $\Hom_c (F,G)$ of maps having bounded image, $G\vec{\otimes} F^*$ with the  space $\Hom_d (F,G)$ of maps having open kernel, and $F^* \otimes^* G$ with $\Hom_f (F,G):= \Hom_c (F,G)\cap \Hom_d (F,G)$. So the vector spaces $\Hom (F,G)$, $\Hom_c (F,G)$, $\Hom_d (F,G)$,  $\Hom_f (F,G)$ carry natural topologies, and (1.1.2) yields  a short exact sequence in $\CT\!op$ $$ 0\to \Hom_f (F,G)\to \Hom_c (F,G)\oplus \Hom_d (F,G) \to \Hom (F,G)\to 0 .\tag 1.2.1$$

Let us describe these topological vector spaces  explicitly. Choose decompositions  $F= F_c \oplus F_d$, $G= G_c \oplus G_d$ where $F_c$, $G_c$ are c-lattices, $F_d$, $G_d$ are d-lattices. Consider the corresponding decomposition of the $\Hom$ spaces into the sum of the four subspaces. One has: 

- $\Hom (F_c ,G_d )=\Hom_c (F_c ,G_d )=\Hom_d (F_c ,G_d )=\Hom_f (F_c ,G_d )$: this is a discrete vector space;

- $\Hom (F_d ,G_c )=\Hom_c (F_d ,G_c ) =\Hom_d (F_d ,G_c )=\Hom_f (F_d ,G_c )$: this is a compact vector space;

- $\Hom (F_c ,G_c )$ $= \Hom_c (F_c ,G_c )$ and $\Hom_d (F_c ,G_c )= \Hom_f (F_c ,G_c )$; the latter topological vector space equals $G_c \otimes F_c^*$ equipped with the ind-topology; 

- $\Hom (F_d ,G_d )= \Hom_d (F_d ,G_d )$ and $\Hom_c (F_d ,G_d )= \Hom_f (F_d ,G_d )$;  the latter topological vector space  equals $F_d^* \otimes G_d $ equipped with the ind-topology. 

{\it A correction to \cite{BD} 2.7.7.} In loc.~cit., certain topologies on the $\Hom$ spaces were considered. The topology on $\Hom (F,G)$ coincides with the above one, but those on $\Hom_c$, $\Hom_d$ and $\Hom_f$ are stronger than the above topologies, and they are   less natural.  Precisely, these topologies differ  at terms $\Hom_d (F_c ,G_c )= \Hom_f (F_c ,G_c )$ and $\Hom_c (F_d ,G_d )= \Hom_f (F_d ,G_d )$ which are considered in \cite{BD} 2.7.7 as discrete vector spaces.
 Exact sequence (1.2) coincides with \cite{BD} (2.7.7.1) (up to signs); it is strongly compatible with either of the topologies. These topologies played an auxiliary role: they were used to define (in \cite{BD} 2.7.8)  the topology on the Tate extension $\fg\fl(F)^\flat$ as the quotient topology for the canonical surjective map $\End_c (F)\oplus \End_d (F) \twoheadrightarrow \fg\fl (F)^\flat$.  Replacing them by the present topologies does not change the quotient topology on $\fg\fl (F)^\flat$. We suggest to discard the topologies on $\Hom_c$, $\Hom_d$, $\Hom_f$ from \cite{BD} 2.7.7, and replace them by those defined  above. 

\medskip

{\bf 1.3. The chiral operations.} As in \cite{BD} 1.1.3, the ! and $*$ tensor structures on $\CT\!op$ can be seen as pseudo-tensor structures with the spaces of operations $P^!_I (\{ V_i \} ,F):= \Hom (\otimes^! V_i ,F)$ and $P^*_I (\{ V_i \} ,F):= \Hom (\otimes^* V_i ,F)$. 

The tensor product $\vec{\otimes}$ is non-commutative, so the definition of the corresponding pseudo-tensor structure requires an appropriate induction:

{\it A  digression: } Let $\Sigma$ be the  linear orders operad (see \cite{BD} 1.1.4).
Suppose we have a $k$-category $\CT$ equipped with an associative tensor product (i.e., a monoidal structure) $\otimes$. For a finite  non-empty collection of objects $\{ V_i \}_{i\in I}$ and a linear order $\tau \in \Sigma_I$, i.e., $\tau : \{ 1,\ldots ,n\} @>{\sim}>> I$, we set $\otimes^\tau V_i := V_{\tau (1)} \otimes \ldots \otimes V_{\tau (n)}$. We
 define the {\it induced} pseudo-tensor structure on $\CT$  by formula $P^{\otimes}_I (\{ V_i \} ,F):= 
\mathop\oplus\limits_{\tau \in \Sigma_I} P^{\otimes}_I (\{ V_i \} ,F)^\tau $ where $P^{\otimes}_I (\{ V_i \} ,F)^\tau := 
\Hom ( \otimes^\tau V_i ,F). $ The composition of  operations is 
defined in the evident way so that the $\Sigma$-grading is compatible with the composition.

If $\otimes$ is actually commutative, then $P^{\otimes}_I (\{ V_i \} ,F)= \Hom (\otimes V_i ,F)\otimes \CA ss_I$ where $\CA ss = k[\Sigma ]$ is the associative algebras operad (see \cite{BD} 1.1.7). The construction is functorial with respect to natural morphisms between $\otimes$'s.

Applying this to $\CT\!op$ and $\vec{\otimes}$, we get the {\it chiral operations}  $P^{ch} := P^{\vec{\otimes}}$; those from the $\tau$-component of $P^{ch}$ are called 
{\it chiral$^\tau$ operations}. They
 define the chiral pseudo-tensor structure $\CT\!op^{ch}$ on $\CT\!op$.
This pseudo-tensor structure is representable (see \cite{BD} 1.1.3) with the pseudo-tensor product
 $\otimes^{ch} V_i =\mathop\oplus\limits_{\tau \in \Sigma_I} \vec{\otimes}^\tau V_i$.

We see that (1.1.1) gives rise to natural transformations $$ P^!_I (\{ V_i \} ,F)\otimes \CA ss_I \to P^{ch}_I (\{ V_i \} ,F)\to P^*_I (\{ V_i \} ,F)\otimes\CA ss_I  \tag 1.3.1$$ compatible with the composition of operations, i.e., the identity functor of $\CT\!op$ lifts naturally to pseudo-tensor functors $\CT\!op^! \otimes\CA ss \to \CT\!op^{ch} \to \CT\!op^* \otimes\CA ss$. 

{\it Remarks.} (i) The arrows in (1.3.1) are injective. If all the $V_i$ are discrete, then they are isomorphisms.

(ii) The endofunctor of $\CT\!op$ which assigns to a topological vector space $V$ same $V$ considered as a discrete space lifts in the obvious manner to a faithful pseudo-tensor endofunctor for all the above pseudo-tensor structures.

Composing (1.3.1) from the left and right with the standard embedding $\CL ie  \hookrightarrow \CA ss$ and  projection $\CA ss \twoheadrightarrow \CC om$,\footnote{Here $\CL ie$ is the Lie algebras operad and $\CC om$ is the commutative algebras operad, i.e, the unit $k$-operad.} we get  natural morphisms (cf.~\cite{BD} 3.2.1, 3.2.2) $$P^!_I (\{ V_i \} ,F)\otimes \CL ie_I \to P^{ch}_I (\{ V_i \} ,F)\to P^*_I (\{ V_i \} ,F). \tag 1.3.2$$

By (1.1.2), in case of two arguments we get an exact sequence $$ 0\to P^!_2 (\{ U,V \} ,F)\otimes\CL ie_2 \to P^{ch}_2 (\{ U,V\} ,F)\to P^*_2 (\{ U,V\} ,F)  \tag 1.3.3$$

{\bf 1.4. Chiral algebras in the topological setting.}
Let $A$ be a topological vector space. 
A  {\it Lie$^{ch}$ algebra} (or {\it non-unital chiral algebra}) structure on $A$ is a Lie bracket $\mu_A : A\otimes^{ch} A\to A$ for the chiral pseudo-tensor structure;  we call such $\mu_A$ a {\it chiral product} on $A$.

Consider the map $P^{ch}_2 (\{ A,A\} ,A)\to \Hom (A\vec{\otimes}A,A)$ which assigns to a binary chiral operation its first component.

\proclaim{\quad Proposition} This map establishes a bijection between the set of chiral products on $A$ and the set of associative  products $A\vec{\otimes}A \to A$.
\endproclaim

{\it Proof.} Our map yields a bijection between the set of skew-symmetric binary chiral operations and $\Hom (A\vec{\otimes} A,A)$. It remains to check that a skew-symmetric chiral operation $\mu$ satisfies the Jacobi identity if and only if the corresponding product $A\vec{\otimes}A\to A$ is associative. 

By Remark (ii) from 1.3  it suffices to consider the case of discrete $A$, so now $A$ is a plain vector space. Let $P(A)$ be the operad of polylinear endooperations on $A$, i.e., $P(A)_I := \Hom (A^{\otimes I},A)$. By Remark (i) from 1.3 the operad $P^{ch}(A)$ of chiral endooperations of $A$ equals $P(A)\otimes \CA ss$. 

Let $\cdot_{as} \in \CA ss_2$ be the standard associative binary operation. We want to prove that a binary operation $\cdot_A \in P(A)_2$ is associative if and only if the binary operation $\cdot_A \otimes\cdot_{as} - \cdot^t_A \otimes \cdot^t_{as}$ in $ P(A) \otimes \CA ss$ is a Lie bracket (here ${}^t$ is the transposition of arguments). Indeed, if $\cdot_A$ is associative, then such is $\cdot_A \otimes\cdot_{as}$, hence the  commutator $\cdot_A \otimes\cdot_{as} - \cdot^t_A \otimes \cdot^t_{as}$ is a Lie bracket. The converse statement is a simple computation left to the reader.   \hfill$\square$

{\it Remark.} Here is a reformulation of the proposition in a  non-commutative algebraic geometry style. We forget about the topologies. Let $A$ be any vector space.  Suppose we have an associative algebra structure on $A$.   Then for any test associative algebra $R$ the tensor product $A\otimes R$ is naturally an associative algebra, hence a Lie algebra. In other words, let $\phi_A$ be the functor $ R\mapsto A\otimes R$ from the category of associative algebras to that of vector spaces; then an associative algebra structure on $A$ yields a Lie algebra structure on $\phi_A$ (i.e., a lifting of $\phi_A$ to the category of Lie algebras). The claim is that this  establishes a bijection between the set of associative algebra structures on $A$ and that of  Lie algebra structures on $\phi_A$.

For $A\in\CT\!op$ a {\it topological associative algebra} structure on $A$ is a continuous associative bilinear product $A\times A \to A$, i.e., an associative product $A\otimes^* A \to A$.
One can sum up the above proposition as follows (the equivalence of (ii), (iii), and (iv) is Remark (ii) in 1.1):

\proclaim{\quad Claim}  For $A\in\CT\!op$ the following structures on $A$ are equivalent:

(i) A non-unital chiral algebra structure $\mu_A$;

(ii) An associative product $\cdot_A :A\vec{\otimes} A \to A$;

(iii) A topological associative algebra structure such that the open left ideals form a base of the topology of $A$. 

(iv) An associative algebra structure such that the open left ideals form a base of the topology of $A$ and the corresponding Lie bracket is continuous.
 \hfill$\square$
\endproclaim

{\it A correction to \cite{BD} 3.6.1:} In loc.~cit.~this was recklessly called ``topological associative algebra" structure; we suggest to rescind this confusing terminology.

{\it A linguistic comment.} The term ``chiral" refers to the breaking of symmetry between
the right and left movers in physics and is rather awkward in the ``purely holomorphic" setting of \cite{BD} (where only one type of movers is present). In the topological setting it looks more suitable for we consider associative products whose left-right asymmetry is enforced by the topology (though now it has to do rather with the time ordrering of physicists).

One often constructs topological chiral algebras using the next corollary:

\proclaim{\quad Corollary} Let $T$ be an associative algebra  equipped with a (not necessary complete) linear topology. Suppose that 

- the open left ideals form a base of the topology;

- for every $r$ from a set of associative algebra generators  the endomorphisms $t\mapsto [r,t] := rt-tr$ of $T$ are continuous. 

Then the completion of $T$ is a chiral algebra.
\endproclaim

{\it Proof.} The endomorphisms $t\mapsto at$ of $T$ are continuous for any $a\in T$ by the first property, so the second property amounts to continuity of the endomorphisms $t\mapsto tr$  for $r$ from our system of generators, hence they are continuous for every $r\in T$. Thus the completion of $T$ carries the structure from (ii) of the claim. \hfill$\square$

\medskip

Suppose we have a topological Lie algebra $L$, i.e., $L$ is a topological vector space equipped with a continuous Lie bracket (i.e., a Lie bracket with respect to $\otimes^*$). We say that $L$ is a {\it topological Lie$^*$ algebra} if it satisfies the following technical condition: 
 the open Lie subalgebras form a base of the topology of $L$. Notice that this condition  holds automatically if $L$ is a Tate vector space.\footnote{Proof: if $P\subset L$ is a c-lattice,  then its normalizer $Q$ is an open Lie subalgebra of $L$, hence $P\cap Q$ is also an open Lie subalgebra.}

Now the second arrow in (1.3.2) transforms any chiral product $\mu_A$ into a continuous Lie bracket $[\enspace ]_A$. Equivalently, $[\enspace ]_A$ is the commutator for the associative product $\cdot_A$.
By (iii) of the claim, $[\enspace ]_A$ is a Lie$^*$ algebra structure on $A$. 

Our  $\mu_A$ is said to be {\it commutative} if $[\enspace ]_A =0$. By by Lemma from 1.1 (and the proposition above) a commutative chiral product amounts to a commutative$^!$ algebra structure, i.e., a commutative and associative product $\cdot_A : A\otimes^! A \to A$.

We say that a Lie$^{ch}$ algebra structure is unital if such is the corresponding associative algebra structure. Such  algebras are referred to simply as {\it topological chiral algebras.} 

For a topological chiral algebra $A$ a {\it discrete $A$-module} is a left unital $A$-module $M$ (we consider $A$ as a mere associative algebra now) such that the action $A\times M\to M$ is continuous (we consider $M$ as a discrete vector space). Equivalently, this is a discrete unital left $A$-module with respect to  $\vec{\otimes}$ monoidal structure. Denote by $A\,$mod$^\delta$ the category of discrete $A$-modules.

Let $\phi : A\,$mod$^\delta \to \CV ec$ be the forgetful functor (which assigns to a discrete $A$-module  its underlying vector space). Then $A$ recovers from $(A\,$mod$^\delta ,\, \phi )$:

\proclaim{\quad Lemma}  $A$  equals the topological associative algebra of endomorphisms of $\phi$. \hfill$\square$
\endproclaim

{\bf 1.5.}
We denote by $\CA ss (\CT\! op^* )$ the category of topological associative unital algebras, 
by $\CC\CA (\CT\!op)$ the category of topological chiral algebras, and by
$\CA ss (\CT\!op^! )$ that of  associative unital algebras with respect to $\otimes^!$. 

As we have seen in 1.4, the above structures on $A\in\CT\!op$ are the same as an associative product $A\otimes A\to A$ that satisfies stronger and stronger continuity conditions. Thus we have fully faithful embeddings $$\CA ss (\CT\!op^! ) \hookrightarrow\CC\CA (\CT\!op)\hookrightarrow \CA ss (\CT\!op^* ).
\tag 1.5.1$$ These embeddings admit left adjoint functors $$ \CA ss (\CT\!op^! ) \leftarrow\CC\CA (\CT\!op)\leftarrow \CA ss (\CT\!op^* ). \tag 1.5.2$$
Namely, for $A\in \CA ss (\CT\!op^* )$ the corresponding chiral algebra $A^{ch}$ is the completion of $A$ with respect to the topology whose base is formed by open left ideals, and for $B\in\CC\CA (\CT\!op)$ the corresponding associative$^!$ algebra is the completion of $B$ with respect to the topology whose base is formed by open two-sided ideals (see Remark (ii) in 1.1).

{\it Remark.} For $B\in \CA ss (\CT\!op^* )$ the category of left unital discrete $A$-modules coincides with $B^{ch}$mod$^\delta$.

The forgetful functor from either of the categories of topological algebras above to $\CT\!op$ which sends a topological algebra to the underlying topological vector space also admits left adjoint. For $V\in \CT\!op$ the corresponding algebras are denoted by $T^* V$, $T^{ch} V$, and $T^! V$. These are completions of the ``abstract" tensor algebra $T V:= \mathop\oplus\limits_{n\ge 0} V^{\otimes n}$ with respect to the topology whose base is formed by all subspaces of type $f^{-1}(U)$, where $f: T V \to A$ is a morphism of associative algebras such that $f|_V : V\to A$ is continuous, $A$ is an algebra of our class, and $U\subset A$ is an open subspace.

{\it Remarks.} (i) Consider T $:= \mathop\oplus\limits_{n\ge 0} V^{\otimes^* n}$. This is a topological vector space (with the inductive limit topology) and an associative algebra, but the product need not be continuous. If $V$ is a Tate space, then the product is continuous if $V$ is either discrete or compact, and not continuous otherwise.\footnote{To check the latter assertion, use the next fact (applied to $P:= V^{\otimes^* 2} $, $T_i := V^{\otimes^* i}$): if $P, T_0 ,T_1,\ldots \in\CT\! op$ are non-discrete and $ P$ is not Tate, then $\oplus\, ( P\otimes^* T_i )\neq P\otimes^* (\oplus T_i )$.}

(ii) One obtains  $T^{ch} V$ and $T^! V$ by applying to $T^* V$ the functors from (1.5.2). Certainly, 
$T^! V=\limleft T (V/P)$ where $P$ runs the set of an open subspaces of $V$ and $T(V/P)$ is its plain tensor algebra (which is a discrete vector space).

{\it A correction to \cite{BD} 3.6.1:} In loc.~cit.~there is a wrong claim that $T^{ch} V$ equals
the direct sum $\vec{T} V := k\oplus V \oplus V^{\vec{\otimes} 2} \oplus \ldots$ (equipped with the direct limit topology).

{\bf 1.6. A topological $\CD$-module setting.} Let $X$ be our curve and $x\in X$ be a point. We will consider pairs $(M, \Xi_M )$ where $M$ is a $\CD$-modules on $X$and $\Xi_M$ is a topology on $M$ at $x$ (see \cite{BD} 2.1.13); such pairs form  a $k$-category $\CM (X,\CT\!op_x )$.  Denote by $\hat{h}_x (M,\Xi_M )$ the completion of the de Rham cohomology stalk $h(M)_x$ with respect to our topology (see loc.~cit.); we get a functor $\hat{h}_x : 
\CM (X,\CT\!op_x )\to \CT\!op$.

We will extend $\hat{h}_x$ to a pseudo-tensor functor with respect to the !, $*$ and chiral polylinear structures. In order to do this, one needs to explain which operations between  $\CD$-modules are continuous with respect to our topologies. 

(a) {\it $*$ operations} (cf.~\cite{BD} 2.2.20). Let $(M_i ,\Xi_{M_i} )$, $i\in I$, be a finite non-empty collection of objects of $\CM (X,\CT\!op_x )$.

\proclaim{\quad Lemma} The $\CD_{X^I}$-module $\boxtimes M_i$ carries a natural topology $\boxtimes^* \Xi_{M_i}$ at $\Delta^{(I)}(x)=(x,\ldots ,x)\in X^I$ such that the corresponding completion of the de Rham cohomology stalk $h (\boxtimes M_i )_{(x,\ldots ,x)}$ equals $\otimes^* \hat{h}_x (M_i ,\Xi_{M_i} )$.
\endproclaim

{\it Proof.} We want to assign in a natural way to every discrete quotient $T$ of $\otimes^* \hat{h}_x (M_i ,\Xi_{M_i} )$ a certain quotient of $\boxtimes M_i$ equal to $i_{(x,\ldots ,x)*} T$. Since $h_x$ commutes with inductive limits, we can assume that each $M_i$ is a finitely generated $\CD_X$-module. Then $
 \hat{h}_x (M_i ,\Xi_{M_i} )$ are all compact, so we can assume that $T= \otimes T_i$ where $T_i$ are discrete (finite-dimensional) quotients of $ \hat{h}_x (M_i ,\Xi_{M_i} )$. So $i_{x*} T_i$ is a quotient of $M_i$, hence $i_{(x,\ldots ,x)*} T =\boxtimes i_{x*} T_i$ is a quotient of $\boxtimes M_i$, and we are done.  \hfill$\square$

\medskip

Let $(N,\Xi_N )$ be another object of $ \CM (X,\CT\!op_x )$. A $*$ operation $\varphi \in P^*_I (\{ M_i \}, N)$ is said to be continuous with respect to our topologies if $\varphi : \boxtimes M_i \to \Delta^{(I)}_* N$ is continuous with respect to the topologies $\boxtimes^* \Xi_{M_i}$ and $\Delta^{(I)}_* \Xi_N$. The composition of continuous operations is continuous, so they form a pseudo-tensor structure $ \CM (X,\CT\!op_x )^* $ on $ \CM (X,\CT\!op_x )$. By construction, $\hat{h}_x$ lifts to a pseudo-tensor functor $$\hat{h}_x :  \CM (X,\CT\!op_x )^* \to \CT\!op^* . \tag 1.6.1$$

(b)  ! {\it  operations.} From now on we will consider a full subcategory $ \CM (U_x ,\CT\!op_x )$ of
 $\CM (X,\CT\!op_x )$ formed by those pairs $(M, \Xi_M )$ that $M = j_{x*} j_x^* M$ where $j_x $ is the embedding $U_x := X\smallsetminus \{ x \} \hookrightarrow X$.  Suppose that our $(M_i ,\Xi_{M_i} )$ lie in this subcategory.
 
 \proclaim{\quad Lemma} The $\CD_X$-module $\otimes^! M_i$ carries a natural topology $\otimes^! \Xi_{M_i}$ at $x$ such that the corresponding completion of $h (\otimes^! M_i )_x$ is equal to $\otimes^! \hat{h}_x(M_i ,\Xi_{M_i})$.
  \endproclaim

{\it Proof.} $\otimes^! \Xi_{M_i}$ is the topology with base formed by submodules $\otimes^! P_i $ of $ \otimes^! M_i$ where $P_i \subset M_i$ are open submodules for $\Xi_{M_i}$.  \hfill$\square$

\medskip 

Set $\otimes^! (M_i ,\Xi_{M_i}):= (\otimes^! M_i ,\otimes^! \Xi_{M_i})$. This tensor product makes  $ \CM (U_x ,\CT\!op_x )$ a tensor category which we denote by $ \CM (U_x ,\CT\!op_x )^!$. Its unit object is
 $j_{x*}\omega_{U_x}$ equipped with the topology formed by the open submodule $\omega_X \subset j_{x*}\omega_{U_x}$. The functor $\hat{h}_x$ lifts naturally to a tensor functor $$\hat{h}_x : \CM (U_x ,\CT\!op_x )^! \to \CT\!op^! .\tag 1.6.2$$

(c) {\it Chiral operations.} Let $(M_i ,\Xi_{M_i} )$ be, as above, some objects of $ \CM (U_x ,\CT\!op_x )$, and $j^{(I)}: U^{(I)}\hra X^I$ be the complement to the diagonal divisor. 

\proclaim{\quad Lemma} The $\CD_{X^I}$-module $j^{(I)}_* j^{(I)*} \boxtimes M_i $ carries a natural topology $\boxtimes^{ch} \Xi_{M_i}$ at $(x,\ldots ,x)\in X^I$ such that the corresponding completion of the de Rham cohomology stalk $h ( j^{(I)}_* j^{(I)*}  \boxtimes M_i )_{(x,\ldots ,x)}$ is equal to $\otimes^{ch} \hat{h}_x (M_i ,\Xi_{M_i} )$.
\endproclaim

{\it Proof.} We proceed by induction by $|I|$. Let us choose for each $i\in I$ a $\Xi_{M_i}$-open submodule $P_i \subset M_i = j_{x*} j_x^* M_i$; set $T_i := h(M_i /P_i )_x$, so $M_i /P_i = i_{x*} T_i$. Set $I_i := I\smallsetminus \{ i\}$. The sequence $0\to \boxtimes P_i \to \boxtimes M_i \to \mathop\oplus\limits_{i\in I} (\mathop\boxtimes\limits_{i'\in I_i} M_{i'})\boxtimes i_{x*} T_i \to 0$ is short exact over  $U^{(I)}$. So we have a short exact sequence $$ 0\to j^{(I)}_* j^{(I)*}\boxtimes P_i \to j^{(I)}_* j^{(I)*}\boxtimes M_i \to \mathop\oplus\limits_{i\in I} (j^{(I_i )}_* j^{(I_i )*}\mathop\boxtimes\limits_{i'\in I_i} M_{i'})\boxtimes i_{x*} T_i \to 0. \tag 1.6.3$$

By the induction assumption, for each $i\in I$ the $\CD_{X^{I_i}}$-module $j^{(I_i )}_* j^{(I_i )*}\mathop\boxtimes\limits_{i'\in I_i} M_{i'}$ carries the   topology $\mathop{\boxtimes^{ch}}\limits_{i'\in I_i} \Xi_{M_{i'}}$
at $(x,\ldots ,x)$ with the completed de Rham cohomology stalk $\mathop{\otimes^{ch}}\limits_{i'\in I_i} \hat{h}_x (M_{i'} ,\Xi_{M_{i'}} )$. Let us equip the $\CD_{X^I}$-module $(j^{(I_i )}_* j^{(I_i )*}\mathop\boxtimes\limits_{i'\in I_i} M_{i'})\boxtimes i_{x*} T_i = i_{X^{I_i}\times x *}( j^{(I_i )}_* j^{(I_i )*}\mathop\boxtimes\limits_{i'\in I_i} M_{i'})\otimes  T_i$ with the ind-topology (cf.~Remark in 1.1), so the  completion of its de Rham cohomology stalk equals 
$\mathop{\otimes^{ch}}\limits_{i'\in I_i} \hat{h}_x (M_{i'} ,\Xi_{M_{i'}} )\vec{\otimes} T_i$. By (1.6.3), the product of these topologies can be seen as a topology  on the quotient module $ j^{(I)}_* j^{(I)*}\boxtimes M_i / j^{(I)}_* j^{(I)*}\boxtimes P_i $ which we denote by $\Xi_{(P_i)}$.

  Now our  topology $\boxtimes^{ch} \Xi_{M_i}$ is formed by all $\CD_{X^I}$-submodules 
$P\subset  j^{(I)}_* j^{(I)*}\boxtimes M_i $ such that $P$ contains  the submodule $j^{(I)}_* j^{(I)*}\boxtimes P_i $ for some choice of $\Xi_{M_i}$-open $P_i \subset M_i$ and the image of $P$ in $ j^{(I)}_* j^{(I)*}\boxtimes M_i / j^{(I)}_* j^{(I)*}\boxtimes P_i $ is $\Xi_{(P_i )}$-open. \hfill$\square$

\medskip

Let $(N,\Xi_N )$ be another object of $ \CM (U_x ,\CT\!op_x )$ and $\varphi \in P^{ch}_I (\{ M_i \}, N)$ be a chiral operation. We say that $\varphi$ is continuous with respect to our topologies if $\varphi : j_*^{(I)} j^{(I)*}\boxtimes M_i \to \Delta^{(I)}_* N$ is continuous with respect to the topologies $\boxtimes^{ch} \Xi_{M_i}$ and $\Delta^{(I)}_* \Xi_N$. The composition of continuous operations is continuous, so they form a pseudo-tensor structure $ \CM (X,\CT\!op_x )^{ch} $ on $ \CM (U_x ,\CT\!op_x )$. By construction, $\hat{h}_x$ lifts to a pseudo-tensor functor $$\hat{h}_x :  \CM (U_x ,\CT\!op_x )^{ch} \to \CT\!op^{ch} . \tag 1.6.4$$

{\bf 1.7. Example.} Let $A$ be a chiral algebra on $U_x$. As in \cite{BD} 3.6.4, we denote by $\Xi_x^{as}$ the topology on $j_{x*} A$ at $x$ whose base is formed by all chiral subalgebras of $j_{x*} A$ that coincide with $A$ on $U_x$. Set $A^{as}_x := \hat{h}_x (j_{x*} A, \Xi_x^{as})$.

\proclaim{\quad Lemma} The chiral product $\mu_A $ is  $\Xi_x^{as}$-continuous.  \hfill$\square$
\endproclaim

Therefore, according to (1.6.4), $A^{as}_x$ is a topological chiral algebra. We leave it to the reader to check that the assocciative product on $A^{as}_x $ coincides with the one defined in \cite{BD} 3.6.6.

\head
2.  Topological cdo.
\endhead

{\bf 2.1. The classical limit: from chiral to coisson algebras.} Let $A$ be any topological vector space. Below a {\it filtration} on $A$ always means an increasing filtration  $A_0 \subset A_1 \subset \ldots$   by closed vector subspaces of $A$ such that  $A_\infty :=\cup A_i$ is dense in $A$. We have a graded topological vector space $\gr_\cdot A$  with components $\gr_i A := A_i /A_{i-1}$.

 The vector space $ \oplus \gr_i A$ carries a natural topology whose base is formed by subspaces $\oplus (P\cap A_i )/(P\cap A_{i-1})$ where $P $ is an open subspace of $A$; we denote by $\widehat{\gr}A$  the completion.

Suppose now that $A$ is a topological chiral algebra and  $A_\cdot$ as above  is a  ring filtration, i.e., $A_i \cdot A_j \subset A_{i+j}$, $1\in A_0$; we call such $A_\cdot$  a {\it chiral algebra filtration}. Then $\widehat{\gr} A$ is naturally a topological chiral algebra. 
 
Our $A_\infty $ is a subring of $A$. Consider a topology on $A_\infty$ formed by all left ideals $I$ in $ A_\infty$ such that  $I\cap A_i$ is open in $A_i$ for every $i$. For any $r\in A_\infty$ the right multiplication endomorphism $a\mapsto ar$ of $A_\infty$ is continuous, so (by Corollary in 1.4) the
 completion $A\tilde{_\infty}$ of $A_\infty$ with respect to this topology is a topological chiral algebra. The embedding $A_\infty \subset A$ extends by continuity to a morphism of topological chiral algebras $A\tilde{_\infty} \to A$. 
 
{\it Definition.} A chiral algebra  filtration $A_\cdot$ is  {\it admissible}  if $A\tilde{_\infty} \iso A$, i.e., a closed left ideal $I\subset A$ is open if (and only if) each $I\cap A_n$ is open in $A_n$.

 {\it Example.} Consider the algebra $k[[t]]$ equipped with the usual topology. Its chiral algebra filtration $k[[t]]_n := k + kt+\ldots + kt^n$ is {\it not}  admissible.

Below  {\it topological coisson algebra} means a  topological vector space $R$ equipped with a Poisson algebra structure such that the Lie bracket  is ($\otimes^*$-) continuous  and the product  is $\otimes^!$-continuous. We also demand  $R$ to be a topological  Lie$^*$ algebra  (see 1.4); equivalently, this means  that open ideals of the commutative algebra $R$ which are Lie subalgebras  form a base of the topology.

Let $A$ be any topological chiral algebra.
A chiral algebra filtration on $A $ is said to be {\it commutative} if $\gr A$ is a commutative algebra. Then  $\widehat{\gr} A$ is a commutative topological chiral algebra, i.e., a commutative$^!$ algebra. The usual Poisson bracket on $\gr A$ extends by continuity to a continuous Lie bracket on $\widehat{\gr}A$, which makes $\widehat{\gr}A$  a topological coisson algebra.

{\bf 2.2. Topological Lie$^*$ algebroids.}
Let $R$ be any (unital) topological commu- tative$^!$ algebra. We denote by $R$mod$^!$ the category of unital $R$-modules in the tensor category $\CT\!op^!$. This is a tensor $k$-category with tensor product ${\mathop\otimes\limits_R}^!$, and an exact category: a short sequence of $R$-modules is exact if it is exact as a sequence in $\CT\!op$ (see 1.1). For $M\in R$mod$^!$ and an open ideal $I\subset R$ we write $M_{R/I}:= (R/I) {\mathop\otimes\limits_R}^! M=M/\overline{IM}$ (here $\overline{IM}$ is the closure  of $IM$); this is a topological $R/I$-module. 

For any $M \in R$mod$^!$ we denote by $\Sym_R^! M$   the universal topological commutative$^!$ $R$-algebra  generated by $M$, i.e.,  $ \Sym^!_R M = \limleft \Sym_{R/J} (M/P)$ where the projective limit is taken with respect all pairs $(J,P)$ where $J\subset R$ is an open ideal and $P\subset M$ an open $R$-submodule such that $J\CL\subset P$. Our $\Sym_R^! M$ carries an evident chiral algebra filtration $(\Sym_R^! M)_a = R\oplus M\oplus  \ldots \oplus \Sym^{!a}_R M$; here $\Sym^{!a}_R M$ is the symmetric power of $M$ in the tensor category $R$mod$^!$.  The filtration  is admissible (due to the universality property of   $\Sym_R^! M$), and $\gr_\cdot 
\Sym_R^! M= 
\Sym^{!\cdot}_R M$ is the universal {\it graded} topological commutative$^!$ $R$-algebra  generated by $M$ in degree 1. 

A {\it topological Lie$^*$ $R$-algebroid} is a topological vector space $\CL$ equipped with a Lie $R$-algebroid structure such that the Lie bracket is ($\otimes^*$-) continuous and  the $R$-action on $\CL$ is $\otimes^!$-continuous. We demand that $\CL$ is a topological Lie$^*$ algebra in the sense of 1.4, or, equivalently, that  open Lie $R$-subalgebroids of $\CL$ form a base of the topology of $\CL$.
  If
$\CL$ is a topological Lie$^*$ $R$-algebroid, then $\Sym_R^! \CL$ is naturally a topological coisson algebra.

{\it Examples.} 
(i) Let $L$ be a topological Lie$^*$ algebra that acts continuously on $R$. Then $L_R := R\otimes^! L$ is naturally a topological Lie$^*$ $R$-algebroid.

(ii)  Let $\Omega_R := \limleft \Omega_{R/I}$ be the topological $R$-module of continuous differentials of $R$. Suppose that $R$ is reasonable, formally smooth, and the topology of $R$ admits a countable base; then $\Omega_R$ is a Tate $R$-module (see \cite{D} Th.~6.2(iii)). Let $\Theta_R$ be the dual Tate $R$-module. Explicitly, $\Theta_R = \limleft \Theta_{R, R/I}$ where for an open ideal $I\subset R$ the topological $R/I$-module $\Theta_{R, R/I}= R/I {\mathop\otimes\limits_R}^! \Theta_R$ consists of all  continuous derivations $\theta : R \to R/I$; the topology of $\Theta_{R, R/I}$ has base formed by $R/I$-submodules that consist of $\theta$ that kill given open ideal $J\subset I$ and  finite subset of $R/J$.
Our $\Theta_R$ is naturally a topological Lie$^*$ $R$-algebroid called  {\it the tangent algebroid} of $R$.

{\bf 2.3. PBW filtrations.}
Let $A$ be  a topological chiral  algebra  equipped with a commutative chiral algebra filtration $A_\cdot$. Set  $R:=A_0$, $\CL:= \gr_1 A =A_1 /A_0$. Then $R$
 is a topological commutative$^!$ algebra and $\CL$ is a topological Lie$^*$ $R$-algebroid. 
 By the universality property we have an evident morphism of  commutative$^!$ $R$-algebras, which is automatically a morphism of topological coisson algebras, 
 $$\Sym^!_R \CL \to \widehat{\gr} A, \tag 2.3.1$$  called the {\it Poincar\'e-Birkhoff-Witt} map. 
 
 \proclaim{\quad Lemma}   Suppose that the filtration    $A_\cdot$ is admissible and each map $\Sym^{!n}_R \CL \to\gr_n A$ is an admissible epimorphism (i.e., an open surjection). 
Then for every open subspace $P\subset A_1$ the closure $\overline{AP}=\overline{A_\infty P}$ of the left ideal  generated by $P$ is open. Such ideals  form a base of the topology of $A$. \endproclaim

{\it Proof.} The second assertion is immediate, once we check the first one. By admissibility,  it suffices to check that   $\overline{A_{n-1}P}$ is open in $A_n$ for each $n$. By induction, we know that $\overline{A_{n-1}P}\cap A_{n-1}$ is open in $A_{n-1}$, i.e., we have an open $V\subset A_n$ such that
$V\cap A_{n-1}\subset \overline{A_{n-1}P}\cap A_{n-1}$. The image $\gr_n \overline{A_{n-1}P}$
of $\overline{A_{n-1}P}$ in $\gr_n A$ is open: indeed,  if an open subspace $T$ of $P$ tends to zero, then $A_{n-1}T$ tends to zero, hence $\gr_n \overline{A_{n-1}P}$ contains the image of $\limleft (\gr_{n-1}A \cdot\gr_1 P )/ (\gr_{n-1}A \cdot\gr_1 T)$    which is an open subspace in $\gr_n A$. Choose $T$ as above such that $A_{n-1}T \subset V$; replacing $V$ by its intersection with the preimage of the open subspace $\gr_n \overline{A_{n-1}T}$ of $\gr_n A$, we can assume that $V= (V\cap A_{n-1}) + (\overline{A_{n-1}T})$, hence $V\subset \overline{A_{n-1}P}$, q.e.d. \hfill$\square$

 {\it Definition.} We say that $A_\cdot$ satisfies {\it the weak PBW property} if $\Sym^{!n}_R \CL \iso \gr_n A$ for each $n$; {\it the strong PBW property} means that the filtration is admissible and (2.3.1) is an isomorphism. Such a filtration is referred to as  {\it weak}, resp.~{\it strong,
 PBW filtration}.

{\it Remarks.} (i) The strong PBW property asserts the existence of large discrete quotients of $A$. Indeed, it amounts to the following two conditions: (a) The filtration is generated by $A_1$, i.e., each $A_n$, $n\ge 1$,  equals the closure of $(A_1 )^n$; (b) Let $P\subset A_1$ be an open subspace such that $\gr_0 P =P\cap R$ is an ideal in $R$ and $\gr_1 P$ is an $R$-submodule of $\gr_1 A$. Then one can find an open ideal $I\subset A$ such that $I\cap A_1 \subset P$ and the projection $\gr_1 A/\gr_1 I \twoheadrightarrow\gr_1 A / \gr_1 P$ lifts to a morphism of  algebras $\gr_\cdot A/\gr_\cdot I \twoheadrightarrow \Sym^{!\cdot}_{R/\gr_0 P} (\gr_1 A /\gr_1 P)$.

(ii) I do not know if every admissible  weak PBW filtration  automatically satisfies the strong PBW property.
  
{\bf 2.4. The chiral envelope of a Lie$^*$ algebra.}
The forgetful functor from the category of topological chiral algebras to that of  Lie$^*$ algebras $(A,\mu_A )\mapsto (A,[\enspace ]_A )$ (see 1.4) admits (as follows easily from 1.5) a left adjoint functor. For a Lie$^*$ algebra $L$ we denote by  $U^{ch} (L)$ the corresponding {\it chiral enveloping} algebra. Explicitly, $U^{ch}(L)$ is the completion of the plain enveloping algebra $U(L)$ with respect to a topology formed by all the left ideals $U(L) P$ where $P\subset L$ is an open vector subspace (it satisfies the conditions of  Corollary in 1.4). 

Our $U^{ch}(L)$ carries a standard filtration $U^{ch}(L)_\cdot$ defined as the completion of the standard (Poincar\'e-Birkhoff-Witt) filtration on $U(L)$. It is admissible and  commutative, so the morphism of Lie$^*$ algebras $L\to \gr_1 U^{ch}(L)$ yields a morphism of topological coisson algebras $$\Sym^! L \to \widehat{\gr}\, U^{ch}(L) . \tag 2.4.1$$

\proclaim{\quad Lemma} (2.4.1) is an isomorphism, i.e.,  $U^{ch}(L)_\cdot$ is a strong PBW filtration.
\endproclaim

{\it Proof.}  By the usual PBW theorem, for every open Lie subalgebra $P\subset L$ one has $\gr_\cdot (U(L)/U(L)P)= \Sym^{!\cdot} (L/P)$ (here the quotient $U(L)/U(L)P$ is equipped with the image of the standard filtration). Such $P$ form a base of the topology of $L$ (by the definition of Lie$^*$ algebra, see 1.4), and we are done.  \hfill$\square$

{\bf 2.5. Chiral extensions of a Lie$^*$ algebroid.}
We want to prove a similar result for a topological Lie$^*$ algebroid. First we need to define its enveloping algebra. This requires (just as in the $\CD$-module setting of \cite{BD} 3.9) an extra structure of chiral extension that we are going to define.

So let $R$ be a topological commutative$^!$ algebra,
 $\CL$ a topological Lie$^*$ $R$-algebroid. 

Consider for a moment $R$ as a commutative algebra and $\CL$ as a Lie $R$-algebroid in the tensor category $\CT\! op^*$. Let $\CL^\flat$ be  a Lie  $R$-algebroid extension of $\CL$ by $R$ in $\CT\! op^*$; below we call such $\CL^\flat$ simply a {\it topological $R$-extension} of $\CL$. Explicitly, our $\CL^\flat$ is an extension of topological vector spaces $$0\to R @>{i}>> \CL^\flat @>{\pi}>>\CL \to 0 \tag 2.5.1$$ together with a Lie $R$-algebroid structure on $\CL^\flat$ such that $\pi$  is a morphism of Lie $R$-algebroids, $i$ is a morphism of $R$-modules, $1^\flat := i(1)$ is a central element of $\CL^\flat$; we also demand that $\CL^\flat$ is a Lie$^*$ algebra and the $R$-action on $\CL^\flat$ is ($\otimes^*$-) continuous. 

{\it Exercise.} $\CL^\flat$ is automatically a topological Lie$^*$ algebra (see 1.4).

Our $\CL^\flat$  is automatically an $R$-bimodule where the right $R$-action is defined by formula $\ell^\flat r =r\ell^\flat +\ell (r)$, where $\ell^\flat \in\CL^\flat $, $r\in R$, $\ell :=\pi (\ell^\flat )$, and $\ell (r)\in R \subset \CL^\flat$; the right $R$-action is continuous as well.

{\it Definition.} (a) $\CL^\flat$ is called a {\it classical $R$-extension} of $\CL$ if $\CL^\flat$ is a topological $R$-algebroid, i.e., the (left) $R$-action on $\CL$ is $\otimes^!$-continuous, so we have $R\otimes^! \CL^\flat \to \CL^\flat$. 
(b) $\CL^\flat$ is called a {\it chiral $R$-extension} of $\CL$ if the left  and  right $R$-actions on $\CL^\flat$ are $\vec{\otimes}$-continuous, i.e., we have $R\vec{\otimes}\CL^\flat \to\CL^\flat$ and $\CL^\flat \vec{\otimes} R \to\CL^\flat$.

{\it Example.} As in 2.3, let $A$ be a chiral algebra equipped with a commutative filtration $A_\cdot$, so $R:= A_0$ is a commutative$^!$ algebra and $\CL:= \gr_1 A$ is Lie$^*$ $A_0$-algebroid. Set $\CL^\flat := A_1$; this is an $R$-extension of $\CL$. The Lie bracket on $A_1$, the left $A_0$-action on $A_1$, and the adjoint action of $A_1$ on $A_0$ make $\CL^\flat$ a topological $R$-extension of the Lie$^*$ $R$-algebroid $\CL$. Since the right $R$-action on $\CL^\flat$ equals  the right $A_0$-action on $A_1$ that comes from the algebra structure on $A$, it is $\vec{\otimes}$-continuous. Therefore $\CL^\flat$ is a chiral $R$-extension of $\CL$.

The topological $R$-extensions of $\CL$ form naturally a Picard groupoid $\CP (\CL )$; the operation is the Baer sum. More precisely, $\CP (\CL )$ is a $k$-vector space in groupoids: For $\CL^{\flat_1}, \CL^{\flat_2}\in\CP (\CL )$ and $a_1 ,a_2 \in k$ the topological $R$-extension $\CL^{a_1 \flat_1 + a_2 \flat_2}$ is defined as the the push-out of $0\to R\times R \to \CL^{\flat_1} \times_{\CL}\CL^{\flat_2} \to\CL \to 0$ by the map $R\times R \to R$, $(r_1 ,r_2 )\mapsto a_1 r_1 + a_2 r_2$; the Lie $R$-algebroid structure on  it  is defined by the condition that the canonical map  $\CL^{\flat_1} \times_{\CL}\CL^{\flat_2} \to \CL^{a_1 \flat_1 + a_2 \flat_2}$ is a morphism of Lie $R$-algebroids.
Let $\CP^{cl} (\CL )$, $\CP^{ch}(\CL )\subset \CP (\CL )$ be the subgroupoids of  classical and chiral $R$-extensions.

\proclaim{\quad Lemma}  $\CP^{cl}(\CL )$ is  a Picard subgroupoid (actually, a  $k$-vector subspace) of $\CP (\CL )$. If $\CP^{ch} (\CL)$ is non-empty, then it is a $\CP^{cl}(\CL )$-torsor. 
\endproclaim

{\it Proof.} Let $\CL^\flat$ be a topological $R$-extension of $\CL$; fix $\lambda \in k$. For $\ell^\flat \in \CL^\flat$, $r\in R$ set $\ell^\flat \cdot_\lambda r =r\ell^\flat + \lambda \ell (r)$. The operation $\cdot_\lambda$ is a right $R$-module structure on $\CL^\flat$ (which commutes with the left $R$-module structure). Notice that $\cdot_1$ is the old right $R$-action on $\CL^\flat$, and $\cdot_0$ is the left $R$-action.

Suppose we have $\CL^{\flat_i}\in \CP (\CL )$ and $\lambda_i ,a_i \in k$; here $i=1,2$.  The right $R$-actions $\cdot_{\lambda_i}$ on $\CL^{\flat_i}$ yield a right $R$-action on $\CL^{\flat_1}\times_{\CL} \CL^{\flat_2}$, hence a right $R$-action on $\CL^{a_1 \flat_1 + a_2 \flat_2}$ such that the canonical map $\CL^{\flat_1}\times_{\CL} \CL^{\flat_2} \to \CL^{a_1 \flat_1 + a_2 \flat_2}$ is a morphism of right $R$-modules. The latter right $R$-action on $\CL^{a_1 \flat_1 + a_2 \flat_2}$ clearly equals  $\cdot_{a_1 \lambda_1 + a_2 \lambda_2}$.

If the  right $R$-actions $\cdot_{\lambda_i}$ on $\CL^{\flat_i}$  are $\vec{\otimes}$-continuous, then the  right $R$-actions on $\CL^{\flat_1}\times_{\CL} \CL^{\flat_2}$ and $\CL^{a_1 \flat_1 + a_2 \flat_2}$ are $\vec{\otimes}$-continuous as well. Therefore the right $R$-action $\cdot_{a_1\lambda_1 + a_2 \lambda_2}$ on $\CL^{a_1 \flat_1 + a_2 \flat_2 }$ is $\vec{\otimes}$-continuous.

Let us call $\CL^\flat $ a {\it $\lambda$-chiral $R$-extension} if the left $R$-action $R\otimes \CL^\flat \to \CL^\flat$ and the right action $\cdot_\lambda : \CL^\flat \otimes R \to \CL^\flat$ are $\vec{\otimes}$-continuous. E.g.,  0-chiral extension is the same as  classical extension, and 1-chiral extension is the same as chiral extension. We have checked that if $\CL^{\flat_i} $ are $\lambda_i$-chiral $R$-extensions,
 then $\CL^{a_2\flat_1 + a_2 \flat_2}$ is an $a_1\lambda_1 + a_2 \lambda_2$-chiral extension. In particular: If $\CL^{\flat_i}$ are classical extensions, then $\CL^{a_1\flat_1 + a_2 \flat_2}$ is a classical extension for every $a_i \in k$. If $\CL^{\flat_1}$ is a classical extension, $\CL^{\flat_2}$ is a chiral one, then $\CL^{\flat_1 +\flat_2}$ is a chiral extension. If $\CL^{\flat_i}$ are chiral extensions, then $\CL^{\flat_1 -\flat_2}$ is a classical extension. We are done.
\hfill$\square$

{\it Exercises.} (i) A classical extension $\CL^\flat$  of $\CL$ 
 is a chiral extension if and only if for every open ideal $I\subset R$ there is an open ideal $J\subset R$ such that $\CL (J)\subset I$. (ii) For a topological $R$-extension $\CL^\flat$ let $\CL^{\flat \tau}$ be the ``inverse" $R$-extension: so we have an identification of topological vector spaces $\CL^{\flat\tau}@>{\sim}>> \CL^\flat$ which commutes with  $i$'s and anticommutes with  $\pi$'s in (2.5.1), interchanges the left and right $R$-module structures, and identifies the Lie bracket with on $\CL^{\flat\tau}$ with minus Lie bracket on $\CL^\flat$. Show that if $\CL^\flat$ is a classical extension, then $\CL^{\flat \tau}$ is a chiral extension if and only if for every open ideal $I\subset R$ there is an open $P\subset \CL$ such that $P(R)\subset I$.

\medskip

{\bf 2.6. The enveloping algebra of a chiral Lie algebroid.}  Let $\CL^\flat$ be a chiral extension of a Lie$^*$ $R$-algebroid $R$; we call a triple $(R,\CL ,\CL^\flat )$ a {\it topological chiral Lie algebroid}. These objects  form naturally a category $\CC\CL (\CT\! op )$.

 Let $\CC\CA^{fc}(\CT\! op )$ be the category of topological chiral algebras equipped with a commutative filtration. By Example in 2.5, we have a functor $\CC\CA^{fc}(\CT\! op )\to \CC\CL (\CT\! op )$,   $(A,A_\cdot )\mapsto (A_0 ,\gr_1 A ,A_1 )$.

 \proclaim{\quad Proposition} This functor admits a left adjoint $\CC\CL (\CT\! op )\to \CC\CA^{fc}(\CT\! op )$.
 \endproclaim 
 
 We denote this adjoint functor by $(R,\CL ,\CL^\flat )\mapsto U^{ch}_R (\CL )^\flat  $ and call $U^{ch}_R (\CL )^\flat$ the {\it chiral enveloping algebra} of $\CL^\flat$. The commutative filtration $U^{ch}_R (\CL )^\flat_\cdot$ on $U^{ch}_R (\CL )^\flat$ is referred to as the {\it standard} filtration.
 
{\it Proof.} Let us define $U_R^{ch}(\CL )^\flat$ as a universal chiral algebra equipped with  a continuous map $\varphi^\flat : \CL^\flat \to U_R^{ch}(\CL )^\flat$ such that $\varphi^\flat$ is a morphism of Lie algebras, its  restriction to $R\subset \CL^\flat$ is a morphism of chiral (or associative) algebras, and $\varphi^\flat$ is a morphism of $R$-bimodules\footnote{It suffices to demand that $\varphi^\flat$ is a morphism of either left or right $R$-modules.} (with respect to $\varphi^\flat |_R $). To construct it explicitly, consider the ``abstract"  enveloping algebra $U_R (\CL )^\flat$  of $\CL^\flat$, i.e.,  take copies of $R$, $\CL$, and $\CL^\flat$ equipped with discrete topologies; then $U_R (\CL )^\flat$ is the corresponding chiral enveloping algebra (its topology is discrete).
 Now $U_R (\CL )^\flat $ carries a linear topology whose base is formed by all left ideals $U_R (\CL )^\flat (P+I)$ where $P\subset L$ and $I\subset R$ are open subspaces. Shrinking $P$, $I$  if necessary, one can assume that $P$ is an open Lie subalgebra of $\CL^\flat$ and $I$ an open $P$-stable ideal of $R$. This topology satisfies
the conditions of Corollary from 1.4, so the corresponding completion is a chiral algebra which equals $U_R^{ch}(\CL )^\flat$ due to the universality property.

The standard  filtration is defined in the usual manner:   $U_R^{ch}(\CL )^\flat_n$ is the closure of the image of $R$ for $n=0$, and  is the closure of the image of $n$th power of the image of
 $\CL^{\flat }$ if $n\ge 1$. It is clearly admissible and commutative. As an object of $\CC\CA^{fc}(\CT\! op )$, our $ U_R^{ch}(\CL )^\flat$ evidently satisfies the universality property of the statement, and we are done.
\hfill$\square$

{\it Remark.} The above explicit construction of $U_R^{ch}(\CL )^\flat$ implies that a
 discrete $U^{ch}_R (\CL )^\flat$-module is the same as a vector space $M$ equipped with a continuous Lie algebra action of $\CL^\flat$ such that the action of $R\subset \CL^\flat$ is  a unital $R$-module structure on $M$, and for $r\in R$, $\ell^\flat \in\CL^\flat$, $m\in M$ one has $(r\ell^\flat )m = r (\ell^\flat m)$.

{\bf 2.7. Rigidified chiral extensions.}
Let $R$ be a topological commutative$^!$ algebra, $L$ a topological Lie$^*$ algebra that acts on $R$ in a continuous way. Then $L_R := R\otimes^! L$ is naturally a topological Lie$^*$ $R$-algebroid.

\proclaim{\quad Lemma} There is a unique, up to a unique isomorphism, chiral extension $L_R^\flat$ equipped with a Lie$^*$ algebra morphism $L \to L_R^\flat$ which lifts the embedding $L \hookrightarrow L_R$. 
\endproclaim
 
{\it Proof.} We  construct $\CL^\flat$ explicitly; the uniqueness is clear from the construction.

 Forgetting for a moment about the topologies, consider the $L$-rigidified Lie $R$-algebroid $L^\delta_R = R\otimes L$ and its trivialized $R$-extension $L_R^{\sharp}=   L^\delta_R \oplus R$. Our $L_R^{\sharp}$ contains $L$ as a Lie subalgebra;  as in 2.5, it is naturally an $R$-bimodule. For open subspaces $P\subset L$ and $I\subset R$ let 
 $\langle P,I \rangle \subset L_R^{\sharp}$ be the vector subspace formed by linear combinations  of all vectors $rp , \ell i,  i'  \in L_R^\sharp$ where $p \in P$,  $i, i' \in I$ and $\ell$, $r$ are arbitrary elements of $L$, $R$. These subspaces form a topology on $L_R^\sharp$; we define $L_R^\flat$ as the corresponding completion.

Notice that the subspaces $\langle P,I\rangle$ with $P\subset L$ an open Lie subalgebra  and $I\subset R$ an open ideal  preserved by the $P$-action form a base of the topology on $L_R^\sharp$. For such $P$, $I$ one has $\langle P,I\rangle = RP + L^\sharp_R I$. One also has
$\langle P,I\rangle \cap R =I$, so $L_R^\flat$ is an extension of $R\otimes^! L$ by $R$. The evident morphism $L\to L_R^\flat$ is continuous. 

For $P$, $I$ as above the subspace $\langle P,I\rangle$ is a Lie subalgebra and a left $R$-submodule (i.e., a Lie $R$-subalgebroid) of $L^\sharp_R$. The Lie bracket and the left and right actions of $R$ on $L^\sharp_R$ are continuous with respect to our topology, hence $L^\flat_R$ is a topological $R$-extension of $L_R$ in the sense of 2.5. Clearly it  is  a chiral $R$-extension, q.e.d. \hfill$\square$

{\it Remark.} Here is another description of $L_R^\flat$ as a mere topological extension. The map $ L\otimes R\oplus R\otimes L \to L_R^\sharp$, $\ell\otimes r+r' \otimes \ell '\mapsto \ell r +r'\ell'$, extends by continuity to a continuous morphism $L\vec{\otimes}R \oplus R\vec{\otimes}L \to L_R^\flat$. It identifies $L_R^\flat$ with the quotient of $ L\vec{\otimes}R \oplus R\vec{\otimes}L \oplus R$ modulo the closed subspace generated by vectors $\ell\otimes r - r\otimes\ell -\ell (r)$. Equivalently, consider the extension $0\to L\otimes^* R \buildrel{(+,-)}\over\longrightarrow L\vec{\otimes}R \oplus R\vec{\otimes}L \buildrel{(+,+)}\over\longrightarrow R\otimes^! L \to 0$ (see (1.1.2)); then $0 \to R \to L_R^\flat \to L_R \to 0$ is its its push-out by the $L$-action map $L\otimes^* R \to R$.

{\bf 2.8. The de Rham-Chevalley chiral Lie$^*$ algebroid.} Let $R$ be a reasonable topological algebra, $\CL$ a Tate $R$-module, $\CL^*$ the dual Tate $R$-module (see \cite{D}).

{\it Exercise.} The $R$-module of all continuous $R$-$n$-linear maps $\CL \times\ldots \times \CL \to R$ identifies naturally with the $n$th tensor power of $\CL^*$ in $R$mod$^!$.

Suppose we have a Lie$^*$  $R$-algebroid structure on $\CL$. Forget for a moment about the topologies and consider $\CL$ as a mere Lie $R$-algebroid. We have the corresponding de Rham-Chevalley complex $\CC\tilde{_R}(\CL )$: this is a commutative DGA whose $n$th term equals $\Hom_R (\Lambda_R^n \CL ,R)$ and the differential is given by the usual formula. By the exercise, the continuous maps form a graded subalgebra $\CC =\CC_R (\CL )$ of $\CC\tilde{_R}(\CL )$ which equals
$\Sym_R^{!\cdot} (\CL^* [-1])$.

\proclaim{\quad Lemma} The de Rham-Chevalley differential preserves the subalgebra $\CC$ and is continuous on it, so $\CC$ is a topological commutative DG algebra.  \hfill$\square$
\endproclaim

{\it Exercise.} A Lie$^*$ $R$-algebroid structure on $\CL$ amounts to a differential on the topological graded algebra $ \Sym_R^{!\cdot} (\CL^* [-1])$.

Our  $\CC$ carries a natural topological DG Lie$^*$  $\CC$-algebroid $\CL_\CC$ (cf.~\cite{BD} 3.9.16).
To construct it, consider $\CL$ as a mere Lie$^*$ algebra. It acts on $\CC$ by transport of structure, and this action extends naturally to an action of the contractible Lie$^*$ DG algebra $\CL_\dagger := \CC one (\id : \CL \to\CL )$: namely, the component $\CL [1]$ acts by the evident $R$-linear derivations of $ \Sym_R^{!\cdot} (\CL^* [-1])$. Thus we have the corresponding DG Lie$^*$ $\CC$-algebroid $\CL_{\dagger\CC}$. By construction, $\CL^{<-1}_{\dagger\CC}=0$ and 
$\CL^{-1}_{\dagger\CC}=R\otimes^! \CL$. Let $\CK$ be the closed DG $\CC$-submodule of   $\CL_{\dagger\CC}$  generated by $\CK^{-1}\subset \CL^{-1}_{\dagger\CC}$ defined as the kernel of the product map $R\otimes^! \CL\to\CL$. Since $\CK^{-1}$ acts trivially on $\CC$ and is normalized by the adjoint action of $\CL_\dagger$, we see that $\CK$ is a DG ideal in the Lie$^*$ $\CC$-algebroid $
\CL_{\dagger\CC}$. The promised $\CL_\CC$ is the quotient $\CL_{\dagger\CC}/\CK$.

Set $\CL_{\dagger \CC+} := \Sym^{!\cdot}(\CL^* [-1])\otimes^!_R \CL [1] $; we consider it as the graded Lie$^*$ $\CC$-algebroid generated by the action of the Lie$^*$ subalgebra $\CL [1]$ of $\CL_{\dagger}$ on $\CC$. The embedding into $\CL_{\dagger \CC}$ identifies it with 
$ \CC\cdot  \CL_{\dagger\CC}^{-1}$, which is a graded submodule of $\CL_{\dagger \CC}$  not preserved by the differential: in fact, the Kodaira-Spencer map $\CL_{\dagger \CC+}\to \CL_{\dagger \CC}/\CL_{\dagger \CC+}$, $\ell \mapsto d(\ell )$mod$\CL_{\CC +}$, is an isomorphism. Similarly, consider $\CK_+ := \CK \cap \CL_{\dagger \CC+}$; then $\CK_+ =  \Sym^{!\cdot}(\CL^* [-1])\otimes^!_R \CK^{-1}$ and 
the Kodaira-Spencer map $\CK_+ \to \CK /\CK_+$ is an isomorphism. Therefore $\CL_{\CC +}$ :=
the closed $\CC$-submodule  of $\CL_\CC$ generated by $\CL_\CC^{-1}$, equals $\Sym_R^{!\cdot}(\CL^* [-1])\otimes^!_R \CL^{-1}_{\CC+} =\Sym_R^{!\cdot}(\CL^* [-1])\otimes^!_R \CL [1]$, and the 
Kodaira-Spencer map $\CL_{\CC + }\to \CL_{\CC}/\CL_{\CC +}$ is an isomorphism. 

{\it Remark.} We see that the Lie$^*$ $R$-algebroid $\CL_\CC^0$ is an extension of $\CL$ by $\CL_{\CC +}^0=\CL^* \otimes^!_R \CL$. It acts naturally on the Tate $R$-module $\CL_\CC^{-1}=\CL$ by the adjoint action. The restriction of this action to the Lie$^*$ $R$-algebra $\CL^* \otimes^!_R \CL$ identifies $\CL^* \otimes^!_R \CL$ with Lie algebra $\fg\fl_R (\CL )$ of continuous $R$-linear endomorphisms of $\CL$ as in 1.2.

\proclaim{\quad Proposition} The topological DG Lie$^*$ $\CC$-algebroid $\CL_\CC$ admits a unique topological chiral DG extension $\CL_\CC^\flat$.
\endproclaim

{\it Proof} (cf.~\cite{BD} 3.9.17). {\it Existence.} Consider the rigidified DG chiral extension $\CL_{\dagger \CC}^\flat$ of the Lie$^*$ $\CC$-algebroid $\CL_{\dagger\CC}$. Let $\tilde{\CK} \subset  \CL_{\dagger \CC}^\flat$ be the  closed DG left $\CC$-submodule of 
$\CL^\flat_{\dagger\CC}$ generated by $\tilde{\CK}^{-1}$ := the preimage of $\CK^{-1}$ by the projection
$\CL_{\dagger \CC}^{\flat\, -1}\to \CL_{\dagger \CC}^{-1}$.(which is an isomorphism since $\CC^{-1}=0$). Then $\tilde{\CK}$ is preserved by the adjoint action of $\CL_\dagger$ (since $\tilde{\CK}^{-1}$ is) and $\tilde{\CK}$ acts trivially on $\CC$, hence $\tilde{\CK}$ is preserved by the adjoint action of $\CL_{\dagger\CC}^\flat$. 

We will show in a moment that $\tilde{\CK}\cap \CC =0$. Then $\CL_\CC^\flat := \CL_{\dagger \CC}^\flat /\tilde{\CK}$ is a chiral DG  extension of $\CL_\CC$, and we are done.

Our   $\tilde{\CK}\cap \CC$    is a DG ideal of $\CC$. It is closed with respect to the action of $\CL [1]\subset \CL_\dagger$, so it suffices to check that $\tilde{\CK}^0 \cap R=0$.  Notice that $\tilde{\CK}^0$ is sum of the closures of $\CC^1 \cdot  \tilde{\CK}^{-1}$ and $d (\tilde{\CK}^{-1})$. The composition $\tilde{\CK}^{-1}\buildrel{d}\over\to \tilde{\CK}^0 \to \CK^0 /\CK^0_+$ is a closed embedding, so $ d (\tilde{\CK}^{-1})$ is a closed subspace of $\CL^{\flat 0}_{\dagger\CC}$ which does not intersect $\CL_{\dagger\CC +}^{\flat 0}$. It remains to prove that the closure of $\CC^1 \cdot  \tilde{\CK}^{-1}\subset   \CL_{\dagger\CC +}^{\flat 0}$ does not intersect $R$. The map $\CC^1 \vec{\otimes}\tilde{\CK}^{-1} \oplus \tilde{\CK}^{-1} \vec{\otimes}\CC^1 \to \CL_{\dagger\CC +}^{\flat 0}$, $c\otimes k + k'\otimes c' \mapsto ck +k'c'$, vanishes on the subspace $\CC^1 \otimes^* \tilde{\CK}^{-1}$ (see (1.1.2)) since $\CK^{-1}$ acts trivially on $\CC$. So, by (1.1.2), the product map $\CC^1 
\vec{\otimes}\tilde{\CK}^{-1}\to \CL_{\dagger\CC +}^{\flat 0}$ extends by continuity to 
 $\CC^1 \otimes^! \tilde{\CK}^{-1}\to  \CL_{\dagger\CC +}^{\flat 0}$ and, by $R$-bilinearity, to
$\mu: \CC^1 \otimes^!_R \tilde{\CK}^{-1}\to  \CL_{\dagger\CC +}^{\flat 0}$. We know that the composition $ \CC^1 \otimes^!_R \CK^{-1}=\CC^1 \otimes^!_R \tilde{\CK}^{-1} \buildrel{\mu}\over\to  \CL_{\dagger\CC +}^{\flat 0}\to \CL_{\dagger\CC +}^{ 0} $ is a closed embedding, so $\mu$ is a closed embedding whose image does not intersect $R$, q.e.d.

{\it Uniqueness.} By the proposition in 2.5, it suffices to show that every classical DG extension $\CL_\CC^{cl}$ of $\CL_\CC$ admits a unique splitting. The uniqueness of the splitting is clear since
$\CL_\CC^{cl\, -1}\iso  \CL_{\CC}^{-1}$ and $\CL_\CC$ is generated by $\CL_{\CC}^{-1}$ as a DG  Lie$^*$ $\CC$-algebroid. To construct one,
 consider  the embedding $\CL [1]\hra \CL_{\CC}^{-1} =\CL_\CC^{cl\, -1}$. It extends to a morphism of complexes $\CL_\dagger =\CC one (\id_\CL )\to  \CL_\CC^{cl}$ which is a morphism of  Lie$^*$ algebras (as an immediate computation shows). By universality, we get a morphism of DG Lie$^*$ $\CC$-algebroids $\CL_{\dagger\CC}\to \CL_\CC^{cl}$ which vanishes on $\CK^{-1}$. Thus it vanishes on $\CK$, hence we get a promised  splitting $\CL_\CC\to \CL_\CC^{cl}$. \hfill$\square$

{\it Remark.} By the $R$-extension property, the Lie$^*$ bracket between $\CL_\CC^{\flat-1}$ and $\CC^1 \subset   \CL_\CC^{\flat 1}$ equals the standard pairing $\CL \times\CL^* \to R$. Therefore the subalgebra of $U_\CC^{ch} (\CL_\CC )^\flat$ generated by these submodules equals the Clifford algebra Cl = Cl$_R (\CL^* \oplus \CL )$ of the Tate $R$-module $\CL^* \oplus \CL$ equipped with the usual hyperbolic quadratic form. It contains the Lie$^*$ $R$-subalgebra $\CL_\CC^{\flat 0}$ which is the {\it Clifford} $R$-extension $\fg\fl_R (\CL )^{Cl}$ of $\fg\fl_R (\CL )$.

{\bf 2.9. $\CD$-modules.}
The most interesting special case of the construction from 2.8 is that of $\CL =\Theta_R$ for $R$ as in Example (ii) in 2.2. Then $\CC $ is the de Rham DG algebra of $R$. 

{\it Remark.} In the vertex (or chiral) algebra setting the enveloping algebra $U^{ch}_\CC (\Theta_R )^\flat$ was first considered in \cite{MSV}
under the name of  the chiral de Rham complex;  for other, essentially equivalent, constructions (in slightly different settings) see \cite{BD}  and \cite{KV}.

 Consider for a moment 
$U^{ch}_\CC (\Theta_R )^\flat$ as a {\it non-graded} topological chiral algebra (which requires an evident completion). The  category of discrete $U^{ch}_\CC (\Theta_R )^\flat$-modules plays the role of the category of $\CD$-modules on the ind-scheme $\Spec R$ (this construction,  mentioned in \cite{D} 6.3.9, generalizes the constructions from \cite{D} and \cite{KV}). 

{\it Exercises.} (i) Suppose that $R$ is discrete. Show that for a discrete $U^{ch}_\CC (\Theta_R )^\flat$-module $M$ the subspace $M^\ell \subset M$ of elements killed by $\Theta_R [1] =\CL_\CC^{\flat \, -1}$ is naturally a {\it left} $\CD$-module on $\Spec \, R$, i.e., an $R$-module equipped with a flat connection. The functor $M\mapsto M^\ell$ is an equivalence between the category of discrete $U^{ch}_\CC (\Theta_R )^\flat$-modules and that of left $\CD$-modules.

(ii) Suppose that $R$ is projective limit of $k$-algebras of finite type, so one has a natural notion of {\it right} $\CD$-module. Show that for a discrete $U^{ch}_\CC (\Theta_R )^\flat$-module $M$ the subspace $M^r \subset M$ of elements killed by $\CC^{\ge 1}$ is naturally a right $\CD$-module, and the functor 
$M\mapsto M^r$ is an equivalence between the category of discrete $U^{ch}_\CC (\Theta_R )^\flat$-modules and that of right $\CD$-modules.

 More generally, suppose $R$ is arbitrary, and we have a fermion module $V$ in the sense of \cite{D} 5.4.1 over the Clifford algebra Cl from the second remark in 2.8. Then $V$ yields a chiral extension $\CL^V$ of $\CL$ defined as follows (cf.~\cite{BD}  3.9.20).
   The Lie$^*$ subalgebra $\CL_\CC^{\flat 0}$ normalizes   Cl $\subset U^{ch}_\CC (\CL_\CC )^\flat$, and its adjoint action on Cl factors through $\CL_\CC^0$.\footnote{In fact, by Remark in 2.8, $\CL_\CC^{0}$ identifies with the Atiyah Lie$^*$ $R$-algebroid of infinitesimal symmetries of the Tate $R$-modules $\Theta_R$ or $\Omega_R$.} Let $\CL_\CC^{0 V}$ be the set of pairs $(\tau,  \tau_V )$ where $\tau \in \CL_\CC^{0}$ and $\tau_V$  is a lifting of $\tau$ to $V$, i.e., a continuous endomorphism of $V$ such that $\tau_V (cv)=\tau (c)v +c\tau_V (v)$ for $c\in$ Cl, $v\in V$. One shows easily that $\CL_\CC^{0 V}$ is naturally a Lie$^*$ $R$-algebroid which is an $R$-extension of $\CL_\CC^0$. Notice that the  action on $V$ of $\fg\fl_R (\CL )^{Cl}\subset$ Cl identifies the restriction of $\CL_\CC^{0 V}$ to $\fg\fl_R (\CL )$ with $\fg\fl_R (\CL )^{Cl}$. Let $ \CL_\CC^{0\flat V}$ be the Baer difference of the $R$-extensions $\CL_\CC^{0 \flat}$ and  $\CL_\CC^{0 V}$; we see that it splits naturally over $\fg\fl_R (\CL )\subset \CL_\CC^0$, so we have a chiral Lie$^*$ $R$-extension $\CL^V :=  \CL_\CC^{0\flat V}/ \fg\fl_R (\CL )$ of $\CL$. 
 
 {\it Exercise.} For a discrete $U^{ch}_\CC (\Theta_R )^\flat$-module $M$ the natural action of 
 $ \CL_\CC^{0\flat V}$ on $M^V :=\Hom_{\text{Cl}}(V,M)$ factors through $\CL^V$, and the functor $M\mapsto M^V$ is an equivalence between the category of discrete $U^{ch}_\CC (\Theta_R )^\flat$-modules and that of discrete modules over the enveloping chiral algebra $U^{ch}_R (\CL )^V$ of $\CL^V$.
 
 {\it Remark.} The key ingredient of $\CO$- and  $\CD$-module theory in the usual finite-dimensional setting is  functoriality of the  derived categories with respect to  morphisms of  varieties.
 The Clifford module picture permits to recover the pull-back functoriality for  appropriately twisted derived categories (the $\Bbb Z$-grading of complexes should be labeled by the dimension $\Bbb Z$-torsor, and the derived category itself should be understood in an appropriate way). This construction is necessary in order to define the notion of $\CO$- or $\CD$-module on orbit spaces such as the moduli space of  (de Rham) local systems on the formal punctured disc. I hope to return to this subject in a joint work with Gaitsgory.

{\bf 2.10. The weak PBW theorem for topological chiral algebroids.} Let $(R,\CL,\CL^\flat )$ be a topological chiral Lie algebroid where $R$ is reasonable.  The morphisms $R\to U^{ch}(\CL )^\flat_0$, $\CL \to \gr_1 U^{ch}(\CL )^\flat$ yield then a morphism of topological coisson algebras (see (2.3.1)) $$ \Sym_R^! \CL \to \widehat{\gr} \,U_R^{ch}(\CL )^\flat .\tag 2.10.1$$

\proclaim{\quad Theorem} If $\CL$ is a flat $R$-module with respect to $\otimes^!$, then $\Sym^{!\cdot}_R \CL \iso \gr_\cdot   U_R^{ch}(\CL )^\flat $.  Thus the standard filtration satisfies the weak PBW property.
\endproclaim

{\it Proof.} We follow \cite{BD} 3.9.13. Set $U^{ch} := U_R^{ch}(\CL )^\flat$ and  $U:= U_R (\CL )^\flat$. 

(a) Suppose our chiral extension admits a rigidification, so we have a Lie$^*$ algebra $L$ acting on $R$, $\CL = L_R =R\otimes^! L$, and $\CL^\flat = L^\flat_R$ is the $L$-rigidified chiral extension. 

Then $U$ is  the enveloping algebra of $(L,R)$, i.e., the quotient of the free associative algebra generated by $L\oplus R$ modulo the relations  saying that the map $L\to U$ is a morphism of Lie algebras and  the map $R\to U$ is a morphism of associative algebras which commutes with the $L$-action (where $L$ acts on $U$ via $L\to U$ and the adjoint action). Our $U^{ch}$ is the corresponding chiral enveloping algebra which can be constructed as follows. Consider the topology on $U$ whose base is formed by all left ideals $U (P+I)$ where $P\subset L$ and $I\subset R$ are open subspaces. Shrinking $P$, $I$  if necessary, one can assume that $P$ is an open Lie subalgebra of $L$ and $I$ an open $P$-stable ideal of $R$. This topology satisfies the conditions from Corollary in 1.4, hence $U^{ch}$ is the completion of $U$ with respect to this topology.

For $(P,I)$ as above the quotient $U/U(P+I)$ coincides with the $L$-module induced from the $P$-module $R/I$. Thus $\gr (U/U(P+I))= \Sym (L/P) \otimes R/I$ by the usual PBW theorem. Passing to the projective limit with respect to $(P,I)$, we see that the standard filtration satisfies the PBW property (actually, the strong one).

(b) Now let $\CL$ be an arbitrary Lie$^*$ $R$-algebroid and $\CL^\flat$ is chiral extension. Let $L$ be a copy of $\CL^\flat$ considered as a mere Lie$^*$ algebra acting on $R$. Consider the corresponding Lie$^*$ $R$-algebroid $L_R = R\otimes^! L$ and its $L$-rigidified chiral extension $L_R^\flat$. We have an evident morphism of Lie$^*$ $R$-algebroids $L_R \twoheadrightarrow \CL$ and its lifting to the chiral extensions $L_R^\flat \twoheadrightarrow \CL^\flat$. The projection $L_R^\flat \to L_R$ identifies the kernels of those morphisms.  
Our $K$ is an $R$-module (in the $\otimes^!$ sense, as a submodule of $L_R$) equipped with a continuous $L$-action (the adjoint one).

Set $\tilde{L}_R := \CC one (K\hra L_R )$, $\tilde{L}^\flat_R := \CC one (K\hra L^\flat_R )$. Then $L_R$ is naturally a DG Lie$^*$ $R$-algebroid, and $\tilde{L}^\flat_R$ is its chiral extension. These structures are uniquely defined by the condition that the  embedding $L_R \to \tilde{L}_R$ is a morphism of Lie$^*$ $R$-algebroids, and that $L^\flat_R \to \tilde{L}^\flat_R$ is a morphism of chiral extensions. Therefore we have a DG chiral algebra $U_R^{ch}(\tilde{L}_R )^\flat$.

Set $\tilde{R} := \Sym^{!\cdot}_R (K[1])$; this is a commutative$^!$ graded topological $R$-algebra whose component in degree $-i$ equals $\Lambda^{! i}_R K$. It carries a natural continuous $L$-action. So we have the ($\Bbb Z$-graded) Lie$^*$ $\tilde{R}$-algebroid $L_{\tilde{R}}$, its $L$-rigidified chiral extension $L_{\tilde{R}}^\flat$, and the corresponding chiral enveloping algebra 
 $U^{ch}_{\tilde{R}} (L_{\tilde{R}})^\flat$.

\proclaim{\quad Lemma} There is a unique isomorphism of $\Bbb Z$-graded topological chiral algebras $$U_R^{ch}(\tilde{L}_R )^\flat \iso U^{ch}_{\tilde{R}} (L_{\tilde{R}})^\flat  \tag 2.10.2$$ which identifies copies of $R$ and $L$ in the degree 0 components, and identifies $K\subset \tilde{L}_R^\flat$ with $K\subset \tilde{R}$ in the components of degree $-1$.
\endproclaim

{\it Proof.} Both  algebras are generated by $R$, $K$, and $L$  with same relations. \hfill$\square$
\medskip

Isomorphism (2.10.2) is identifies the standard filtrations up to a shift by the grading:  one has 
$U^{ch}_R (\tilde{L}_R )^{\flat a}_n \iso U^{ch}_{\tilde{R}} (L_{\tilde{R}})^{\flat a}_{a+n}$. By  (a),  $U_{\tilde{R}}^{ch} (L_{\tilde{R}})^\flat$ satisfies the PBW property. So we have an isomorphism  of graded DG chiral algebras
$$\Sym_R^{!\cdot }\tilde{L}_R \iso \gr_\cdot U_R^{ch}(\tilde{L}_R )^{\flat}. \tag 2.10.3$$

(c) Suppose  that $\CL$ is $R$-flat with respect to $\otimes^!$. Then the projection $\Sym_R^{!\cdot }\tilde{L}_R \to \Sym^{!\cdot}_R \CL$ is a quasi-isomorphism in  the exact category $R$mod$^!$. Therefore
the differential $d$ on $U_R^{ch}(\tilde{L}_R )^{\flat}$ is strictly compatible with the standard filtration. Let  $I_{n} \subset  U_R^{ch}(\tilde{L}_R )^{\flat }_{n}$ be the closure of 
$K U_R^{ch}(\tilde{L}_R )^{\flat}_{n-1} $; then $I_{n} = I_{m}\cap U^{ch}(\tilde{L}_R )^{\flat }_{n}$ for any $m\ge n$, hence  $\gr_\cdot  U^{ch}_R (\CL )^\flat  \buildrel{\sim}\over\leftarrow \gr_\cdot U_R^{ch}(\tilde{L}_R )^{\flat } /\gr_\cdot I \iso \Sym_R^{!\cdot }\CL$, q.e.d.  \hfill$\square$

{\it Exercise.} Suppose the topology of $R$ has a base formed by  open reasonable ideals $I$ that satisfy the next property:  The open Lie$^*$ $R$-subalgebroids $M\subset \CL$ such that $M(I)\subset I$, $M\supset I \CL$, and $\CL/M$ is a flat $R/I$-module, form a base of the topology of $\CL/\overline{I \CL}$. Then the standard filtration on $U_R^{ch}(\CL )^\flat $ is a strong PBW filtration.

\bigskip

\head 
References
\endhead

\ref\key BD
\by\quad  A.~Beilinson and V.~Drinfeld
\book Chiral Algebras
\publ AMS
\publaddr Providence, RI
\yr 2004
\endref

\ref\key D  \by V.~Drinfeld \paper  Infinite-dimensional vector bundles in algebraic geometry: an introduction  \inbook The unity of mathematics \pages  263--304
\bookinfo Progr.~Math. \vol 244 \publ Birkh\"auser \publaddr Boston, MA \yr 2006 \endref
\ref\key G \by  A.~Grothendieck \paper Produits tensoriels topologiques et espaces nucl\' eaires. \jour Mem.~AMS \yr 1955  \issue 16 \endref

\ref\key MSV \by \qquad F.~Malikov, V.~Schechtman, and A.~Vaintrob \paper Chiral de Rham complex \jour Commun.~Math.~Phys. \vol 204 \yr 1999 \pages 439--473\endref

\ref\key KV \by\quad M.~Kapranov and E.~Vasserot \paper Vertex algebras and the formal loop space \jour  Publ.~Math.~IHES \vol 100  \yr 2004 \pages 209--269 \endref

\ref\key S \by J.-P.~
Schneiders \paper Quasi-abelian categories and sheaves \jour M\' em.~Soc.~Math.~Fr. \yr 1999 \issue 76 \endref
\enddocument